\documentclass[12pt]{article}
\usepackage{mathrsfs}
\usepackage{latexsym,amsmath,amssymb,amsfonts,epsfig,graphicx,cite,psfrag}
\usepackage{eepic,color,colordvi,amscd,amsthm}

\newtheorem{thm}{Theorem}[section]

\newtheorem{prop}[thm]{Proposition}

\numberwithin{equation}{section}

\def\pf{\noindent{\it Proof.} }

\def\qed{\nopagebreak\hfill{\rule{4pt}{7pt}}\medbreak}

\begin{document}
\begin{center}
{\large\bf Identities Derived from
Noncrossing Partitions of Type $B$}\\
\end{center}
\begin{center}
William Y. C. Chen$^{1}$, Andrew Y. Z. Wang$^2$ and
Alina F. Y. Zhao$^{3}$\\[6pt]
Center for Combinatorics, LPMC-TJKLC\\
Nankai University, Tianjin 300071, P. R. China\\[5pt]
$^{1}${chen@nankai.edu.cn}, $^{2}${yezhouwang@mail.nankai.edu.cn},
$^{3}${zfeiyan@mail.nankai.edu.cn}
\end{center}

\vskip 3mm \noindent \textbf{Abstract.} Based on weighted
noncrossing partitions of type $B$, we obtain type $B$ analogues of
Coker's identities on the Narayana polynomials. A parity reversing
involution is given for the alternating sum of Narayana numbers of
type $B$. Moreover, we  find  type $B$ analogues of the refinements
of Coker's identities due to Chen, Deutsch and Elizalde. By
combinatorial constructions, we provide type $B$ analogues of three
identities of Mansour and Sun also on  the Narayana polynomials.



\section{Introduction}

The objective of this paper is to give type $B$ analogues of
combinatorial identities on the Narayana polynomials
\cite{Bon,Chen3, Coker,Sul,Sul04}
\[N_n(x)=\sum\limits_{k=1}^{n} N_{n,k} x^k, ~ n\geq1,
\]
where
\[ N_{n,k}= \frac{1}{n}{n\choose k-1}
{n\choose k} \] are called the Narayana numbers.   Let \[
C_n={1\over n+1}{2n\choose n}\] be the $n$-th Catalan number. Using
generating functions, Coker \cite{Coker, Chen3} has derived the
following identities
\begin{eqnarray} \label{c1}
\sum\limits_{k=0}^{n-1}\frac{1}{n}{n\choose k}{n\choose k+1}x^k=
\sum\limits_{k=0}^{\lfloor{(n-1)/2}\rfloor}C_k{n-1\choose 2k}x^k
(1+x)^{n-2k-1},
\end{eqnarray}
\begin{eqnarray}\label{c2}
\sum\limits_{k=0}^{n-1}\frac{1}{n}{n\choose k}{n\choose k+1}x^{2k}
(1+x)^{2(n-1-k)}= \sum\limits_{k=0}^{n-1}C_{k+1}{n-1\choose k}x^k
(1+x)^k.
\end{eqnarray}
Chen, Yan and Yang \cite{Chen3} have given combinatorial
interpretations of these two identities based on weighted  Dyck
paths and $2$-Motzkin paths in answer to a question raised by Coker.
It should be noticed that the identity \eqref{c1} can also be
derived from the following identity due to Simion and Ullman
\cite{Simon}, see also Chen, Deng and Du \cite{Chen0}:
\[\frac{1}{n}{n\choose k}{n\choose k+1}=\sum_{r=0}^{k-1}{n-1 \choose
2r}{n-2r-1 \choose k-1-r}C_r.\]

 Recently, Mansour and Sun \cite{Man} have established the following
three identities on the Narayana polynomials and have given both
algebraic and combinatorial proofs:
\begin{eqnarray}\label{m1}
x^{\frac{n}{2}+1}C_{\frac{n}{2}}& = &
\sum\limits_{k=0}^{n}(-1)^{n-k} {n\choose
k}N_{k+1}(x)(1+x)^{n-k},\\[6pt]
\label{m2} x^{n+2}C_{n+1}& = & \sum\limits_{k=0}^{n}(-1)^{n-k}
{n\choose k}N_{k+1}(x^2)(1-x)^{2(n-k)}, \\[6pt]
\label{m3} C_n & = &
\sum\limits_{k=0}^{n}\frac{2k+1}{2n+1}{2n+1\choose
n-k}N_k(x)(1-x)^{n-k},
\end{eqnarray}
where $C_{\frac{n}{2}}$ is treated as zero if $n$ is odd.

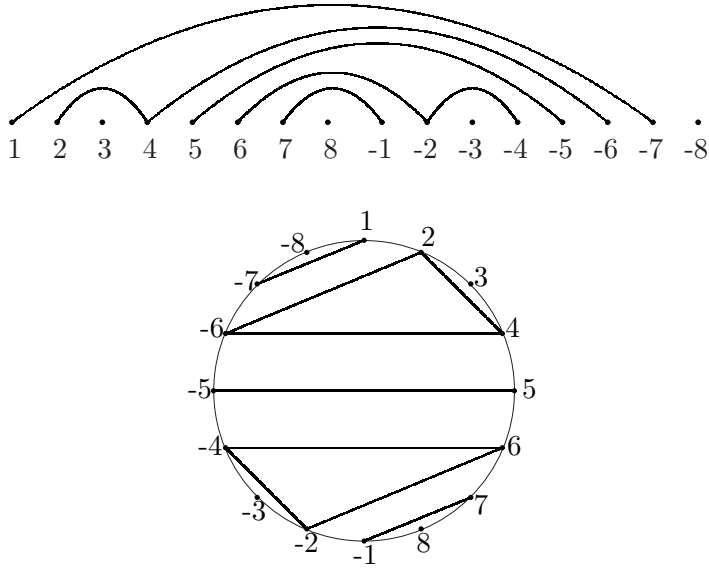
\begin{figure}[h,t]
\setlength{\unitlength}{0.6mm}
\begin{center}
\begin{picture}(160,35)

\put(1,0){\small 1}     \put(11,0){\small 2}    \put(21,0){\small 3}
\put(31,0){\small 4}    \put(41,0){\small 5}    \put(51,0){\small 6}
\put(61,0){\small 7}    \put(71,0){\small 8} \put(81,0){\small{-1}}
\put(91,0){\small{-2}}  \put(101,0){\small{-3}}
\put(111,0){\small{-4}} \put(121,0){\small{-5}}
\put(131,0){\small{-6}} \put(141,0){\small{-7}}
\put(151,0){\small{-8}}

\multiput(2,8)(10,0){8}{\circle*{1}}
\multiput(84,8)(10,0){8}{\circle*{1}}

\qbezier[500](2,8)(73,60)(144,8)  \qbezier[500](12,8)(22,23)(32,8)
\qbezier[500](32,8)(83,50)(134,8) \qbezier[500](42,8)(83,43)(124,8)
\qbezier[500](52,8)(73,30)(94,8)  \qbezier[500](94,8)(104,23)(114,8)
\qbezier[500](62,8)(73,23)(84,8)

\end{picture}
\end{center}

\begin{center}
\setlength{\unitlength}{2mm}
\begin{picture}(20,23)
\put(10,10){\circle{20}}

\put(10,20){\circle*{0.28}}      \put(13.8,19.2){\circle*{0.28}}
\put(17.1,17.1){\circle*{0.28}}  \put(19.2,13.8){\circle*{0.28}}
\put(20.0,10.0){\circle*{0.28}}  \put(19.2,6.2){\circle*{0.28}}
\put(17.1,2.9){\circle*{0.28}}   \put(13.8,0.8){\circle*{0.28}}
\put(10.0,0.0){\circle*{0.28}}   \put(6.2,0.8){\circle*{0.28}}
\put(2.9,2.9){\circle*{0.28}}   \put(0.8,6.2){\circle*{0.28}}
\put(0.0,10.0){\circle*{0.28}}    \put(0.8,13.8){\circle*{0.28}}
\put(2.9,17.1){\circle*{0.28}}    \put(6.2,19.2){\circle*{0.28}}

\qbezier[500](10,20)(6.45,18.55)(2.9,17.1)
\qbezier[500](13.8,19.2)(7.3,16.5)(0.8,13.8)
\qbezier[500](13.8,19.2)(16.5,16.5)(19.2,13.8)
\qbezier[500](19.2,13.8)(10.0,13.8)(0.8,13.8)
\qbezier[500](20.0,10.0)(10.0,10.0)(0.0,10.0)
\qbezier[500](19.2,6.2)(12.7,3.5)(6.2,0.8)
\qbezier[500](6.2,0.8)(3.5,3.5)(0.8,6.2)
\qbezier[500](0.8,6.2)(10.0,6.2)(19.2,6.2)
\qbezier[500](17.1,2.9)(13.55,1.45)(10.0,0.0)

\put(9.7,20.7){\small{1}}        \put(13.8,19.6){\small{2}}
\put(17.3,17){\small{3}}         \put(19.4,13.6){\small{4}}
\put(20.5,9.5){\small{5}}        \put(19.5,5.7){\small{6}}
\put(17.3,1.8){\small{7}}        \put(13.5,-0.7){\small{8}}
\put(9.3,-1.5){\small{-1}}       \put(5.4,-0.7){\small{-2}}
\put(1.9,1.4){\small{-3}}        \put(-1.0,5.7){\small{-4}}
\put(-1.7,9.5){\small{-5}}       \put(-0.9,13.6){\small{-6}}
\put(1.4,16.7){\small{-7}}       \put(4.5,19.2){\small{-8}}

\end{picture}
\end{center}
\caption{The linear and cyclic representations of a type $B$
noncrossing partition.}\label{fig1}
\end{figure}

We obtain type $B$ analogues of the above identities of Coker, and
Mansour and Sun, based on the structure of type $B$ noncrossing
partitions. Recall that a type $B$ partition of $[n]$ is a partition
of the set $\{1,2,\ldots,n,-1,-2,$ $\ldots,-n\}$, which may have a
unique block, it it exists, called the zero block in which $i$ and
$-i$ appear in pairs, such that for any block $B$ of $\pi$, the set
$-B$, obtained by negating the elements of $B$, is also a block of
$\pi$, see, for example, \cite{Arm,Athana,Reiner,Simon4}. Evidently,
the zero block is a union of antipodal pairs $\{i,-i\}$ if it
exists. Moreover, there does not exist any other block $B$ such that
$B=-B$. A type $B$ partition $\pi$ can be represented by a diagram,
with the elements $1,2,\ldots,n,-1,-2,\ldots,-n$ drawn on a
horizontal line in the following order
\begin{equation}\label{order}
 1<2<\cdots < n < -1 <
-2 < \cdots < -n. \end{equation} Accordingly, one can list the
elements of a block in the above order. Suppose that $B=\{i_1, i_2,
\ldots, i_k\}$ is a nonzero block of a type $B$ partition $\pi$. One
may represent this block by a path from $i_1$ to $i_k$ with arcs
drawn above the horizontal line from $i_1$ to $i_2$,  $i_2$ to
$i_3$, and so on. A block with one element is called a singleton
block.  Such a diagram is called the linear representation of a type
$B$ partition. A type $B$ partition is said to be noncrossing if its
diagram contains no crossing arcs, see, Athanasiadis \cite{Athana}.
It is worth mentioning that we may also place the elements $1, 2,
\ldots, n, -1, -2, \ldots, -n$ on a circle, and use a cycle to
represent a block. This is called the cyclic representation of a
type $B$ partition. This leads to  an equivalent definition of
noncrossing partitions of type $B$, see, Reiner \cite{Reiner}. An
illustration of these two representations of a type $B$ noncrossing
partition is given in Figure \ref{fig1}, where \[
\pi=\{1,-7\}\{7,-1\}\{2,4,-6\}\{6,-2,-4\}\{3\}\{-3\}\{5,-5\}\{8\}\{-8\}.\]

In this paper, we shall adopt the linear representation of type $B$
noncrossing partitions. The  set of type $B$ noncrossing partitions
on $[n]$ will be denoted by $NC^B(n)$. It is well known that the
cardinality of $NC^B(n)$ equals ${2n \choose n}$, see, for example,
Reiner \cite[Proposition 6]{Reiner}. Furthermore,  the set of type
$B$ noncrossing partitions of $[n]$ having $k$ pairs of nonzero
blocks will be denoted by $NC^B(n,k)$. The cardinality of
$NC^B(n,k)$, which is known to be ${n\choose k}^2$, has been
referred to as the Narayana number of type $B$ by Fomin and Reading
\cite{Fom}. It also equals the rank-size of the lattice of the
noncrossing partitions of type $B$ of $[n]$ with the refinement
order \cite{Simon4}.  The polynomials
\[
P_n(x)=\sum\limits_{k=0}^n {n \choose k}^2x^k,~n\geq1
\]
will be   called the Narayana polynomials of type $B$.

This paper is organized as follows.  In Section 2, we give type $B$
analogues of Coker's identities and combinatorial proofs in terms of
weighted type $B$ noncrossing partitions. We find an involution for
the case of the alternating sum of the Narayana numbers of type $B$.
We also provide refinements of Coker's identities of type $B$.
Section 3 is devoted to  type $B$ analogues of three identities due
to  Mansour and Sun \cite{Man}.

\section{ Type $B$ Analogues of Coker's Identities }

This section is concerned with  type $B$ analogues of Coker's
identities. We shall use weighted type $B$ noncrossing partitions as
the underlying combinatorial structure. Roughly speaking, we shall
assign a weight to each point, and the weight of a partition is the
product of the weights of all the points. To be precise, we shall
assign weights only to half of the elements appearing in the
canonical representation of a type $B$ noncrossing partition.  The
following two identities can be regarded as type $B$ analogues of
Coker's identities.

\begin{thm} \label{coker}
For $n\geq 0$,
\begin{eqnarray}\label{i1}
\sum\limits_{k=0}^{n}{n\choose k}^2x^k= \sum\limits_{k=0}^{\lfloor
\frac{n}{2}\rfloor}{n\choose 2k}{2k\choose k}x^k(1+x)^{n-2k},
\end{eqnarray}
\begin{eqnarray}\label{i2}
\sum\limits_{k=0}^{n}{n\choose k}^2x^{2k}(1+x)^{2(n-k)}=
\sum\limits_{k=0}^{n}{n\choose k}{2k\choose k}x^k(1+x)^k.
\end{eqnarray}
\end{thm}

Note that identity \eqref{i1} was first derived by Riordan
\cite{Rio} by using generating functions. Before presenting the
combinatorial proof of the above theorem, we recall a basic property
of the linear representation of a type $B$ noncrossing partition, as
 observed by Athanasiadis \cite{Athana}. By a pure block we mean a block
 that contains only positive elements or only negative elements.
 A block is called a mixed block if it is not pure.

\begin{prop}\label{propair}
Let $\pi$ be a type $B$ noncrossing partition.  In the linear
representation of $\pi$, for any pair of antipodal blocks  $B$ and
$-B$, either one lies entirely on the left of the other, or one is
completely covered (or nested)  by an arc of the other.
\end{prop}

\pf When $B$ is a pure block, the assertion is obvious. We now
assume that $B$ is a mixed block.  Let $i$ be the maximum positive
element  and  $-j$ be the minimum negative element of $B$ according
to the order (\ref{order}). So  $(i,-j)$ is an arc in the linear
representation of $\pi$. If $j>i$, it is easily seen that $-B$ is
nested by the arc $(i,-j)$. Similarly, in the case $j<i$, $B$ is
nested by an arc of $-B$. This completes the proof.  \qed

In light of the above property, for a pair of antipodal blocks $B$
and $-B$, we need only one of them to represent this pair. We shall
choose $B$  such that either it is on the left of $-B$ or it is
nested by an arc of $-B$. Moreover, we shall list the representative
blocks $B_1, B_2, \ldots,B_k$ in increasing order of their minimum
elements. In particular, we  use $B_0$ to denote the set of positive
elements of the zero block. We shall call $B_0/ B_1/ B_2/ \cdots/
B_k$ the canonical representation of $\pi$. It is clear that the
elements appearing in the canonical representation of a type $B$
partition of $[n]$ form a set in which either $i$ or $-i$ appears,
but not both, for any $i \in [n]$. As an example, the representative
blocks of the noncrossing partition $\pi$ in Figure \ref{fig1} can
be read as $B_0=\{5\}, B_1=\{3\}, B_2=\{6,-2,-4\}, B_3=\{7,-1\}$ and
$B_4=\{8\}$.

 From now on,
 we shall use the above canonical
representation $\pi=B_0/ B_1/ \cdots/B_k$ for a type $B$ noncrossing
partition. The elements appearing in the canonical representation,
namely, the elements of $B_0, B_1, \ldots, B_k$,  will be classified
into five types.  Let $i\in B_j$. Then we say that
\begin{enumerate}
  \item $i$ is a \emph{zero point} if  $i\in B_0$;
  \item $i$ is a \emph{singleton} if $B_j$ is a nonzero singleton block, that is, $|B_j|=1$;
  \item $i$ is a \emph{transient point} if $i$ is neither the smallest nor the
        largest element of $B_j$, according to the order (\ref{order});
  \item $i$ is a \emph{departure point} if $i$ is the smallest element of $B_j$ and $|B_j|>1$;
  \item $i$ is a \emph{destination point} if $i$ is the largest element of $B_j$ and
$|B_j|>1$.
\end{enumerate}

For example, suppose that $\pi$ is the partition in Figure
\ref{fig1}. Then $5$ is a zero point; $3$ and $8$ are singletons;
$-2$ is a transient point; $6$ and $7$ are departure points; $-1$
and $-4$ are destination points. Before proving Theorem \ref{coker},
we first present a formula that will be used later.

\begin{prop} \label{pro}
The number of partitions in $NC^B(n)$ with exactly $k$ pairs of
nonzero blocks but no singletons equals to
\[{n\choose 2k}{2k\choose k}.\]
\end{prop}

\pf From the correspondence by Reiner \cite{Reiner} between  type
$B$ noncrossing partitions and  pairs $(L,R)$ of $k$-subsets of
$[n]$, it is not hard to see  that a point $i$ is a singleton of
$\pi$ if and only if $i$ appears in  both $L$ and $R$. Thus a type
$B$ noncrossing partition $\pi$ without singletons is uniquely
determined by a pair $(L,R)$ of disjoint subsets of $[n]$ with equal
cardinality, namely $L,R\subseteq [n]$, $L\cap R=\emptyset$ and
$|L|=|R|$. Moreover, the number of pairs of nonzero blocks equals
the cardinality of $|L|$ and $|R|$. Clearly, there are ${n\choose
2k}{2k\choose k}$ ways to choose $(L,R)$. This completes the proof.
\qed

To make this  paper self-contained, we give a more detailed
description of the procedure to generate a type $B$ noncrossing
partition $\pi$ without singletons from a pair $(L, R)$ of disjoint
$k$-subsets of $[n]$  by using the beautiful construction of Reiner
\cite{Reiner}. First, put $2n$ points on a horizontal line with
labels $1,2,\ldots,n,-1,-2,\ldots,-n$ from left to right in
accordance with the order (\ref{order}). If $l_i\in L$, we replace
the points  $l_i$ and $-l_i$ each by a left parenthesis; if $r_i\in
R$, replace the points $r_i$ and $-r_i$ each by a right parenthesis.
Thus at the positions of the elements $1, 2, \ldots,
n,-1,-2,\ldots,-n$, the pair $(L, R)$ corresponds to a sequences of
$2k$ left parentheses, $2k$ right parentheses and $2n-4k$ points.

It is now important to recall a property of the $2k$ parentheses in
the positions $1, 2, \ldots n$, as discovered by Greene and Kleitman
\cite{Gree} in the construction of the symmetric chain decomposition
of the posets of subsets of a finite set. To be more precise, any
sequence of left parentheses and right parentheses consists of
well-parenthesized segments separated by parentheses which can be
read from left to right as $)\cdots)\,(\cdots($. This is to say that
the unpaired right parentheses appear before the unpaired left
parentheses. Moreover, there is no left or right parentheses between
any well-paired parentheses. Thus any sequence of parentheses can be
decomposed into well-parenthesized segments, separated by a sequence
of right parentheses followed by left parentheses. For example, the
sequence $)\,(\,(\,)\,)\,)\,(\,(\,)\,($ has the following
decomposition into well-parenthesized segments.

\begin{figure}[h,t]
\setlength{\unitlength}{0.8mm}
\begin{center}
\begin{picture}(80,8)

\put(0,5){$)$} \put(10,5){$($} \put(20,5){$($} \put(30,5){$)$}
\put(40,5){$)$} \put(50,5){$)$} \put(60,5){$($} \put(70,5){$($}
\put(80,5){$)$}\put(90,5){$($}

 \put(12,-6){\line(0,1){8}}\put(41,-6){\line(0,1){8}}
\put(12,-6){\line(1,0){29}}
\put(22,-3){\line(0,1){4}}\put(31,-3){\line(0,1){4}}
\put(22,-3){\line(1,0){9}}
\put(71,-3){\line(0,1){4}}\put(81,-3){\line(0,1){4}}
\put(71,-3){\line(1,0){10}}
\end{picture}
\end{center}
\end{figure}

For the well-paired segments at the positive positions, or at the
negative positions, we can easily construct pure blocks. For a pair
of  the form $( \, \cdot \cdot \cdot )$, that is, a pair of
parentheses for which there is no parentheses between them,  we
simply form a pure block by selecting the corresponding elements of
the parentheses and the points between them. After such a block is
formed, we may delete the underlying elements and continue the above
procedure until all well-paired parentheses at the positive
positions or at the negative positions are processed.

 Upon the deletion of  the elements of all pure blocks,  the remaining
parentheses have the following form
\[\underbrace{\cdots)\cdots)\cdots)\cdots(\cdots(\cdots(\cdots}_{+}
|\underbrace{\cdots)\cdots)\cdots)\cdots(\cdots(\cdots(\cdots}_{-}.\]
The positive left parentheses and negative right parentheses can be
well-paired, which will lead to mixed blocks. It can be shown  that
a mixed block $B$ formed in the above procedure must be nested by
its antipodal block. If $(i,-j)$ is a paired parentheses which
yields a mixed block $B$, then $j$ is the largest positive element
of $-B$ and $-i$ is the smallest negative element of $-B$. Clearly,
the mixed block $B$ is nested by the arc $(j,-i)$ of $-B$ since
$j<i$. Moreover, one readily sees that the block $B$ forms a
consecutive segment with respect to the order (\ref{order}) after
removing the pure blocks, whereas the block $-B$ occupies two
consecutive  segments at both ends. This implies that one can
continue this process to get noncrossing blocks. When all the pure
and mixed blocks are obtained, if there are still some elements
left, we collect them together to form the zero block.

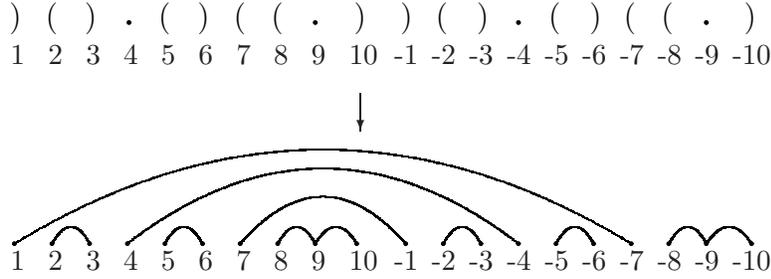
\begin{figure}[h,t]
\setlength{\unitlength}{0.5mm}

\begin{center}
\begin{picture}(200,35)

\put(1,8){\small 1}    \put(11,8){\small 2} \put(21,8){\small 3}
\put(31,8){\small 4}   \put(41,8){\small 5}  \put(51,8){\small 6}
\put(61,8){\small 7}   \put(71,8){\small 8} \put(81,8){\small 9}
\put(91,8){\small 10}

\put(103,8){\small{-1}} \put(113,8){\small{-2}}
\put(123,8){\small{-3}} \put(133,8){\small{-4}}
\put(143,8){\small{-5}} \put(153,8){\small{-6}}
\put(163,8){\small{-7}} \put(173,8){\small{-8}}
\put(183,8){\small{-9}} \put(193,8){\small{-10}}

\put(32.5,19){\circle*{1}} \put(82,19){\circle*{1}}
\put(136,19){\circle*{1}}  \put(186,19){\circle*{1}}

\put(0.7,18){$)$}  \put(10.3,18){$($} \put(20.7,18){$)$}
\put(40.3,18){$($} \put(50.7,18){$)$} \put(60.3,18){$($}
\put(70.3,18){$($} \put(92.3,18){$)$}

\put(104.7,18){$)$}  \put(114,18){$($} \put(125,18){$)$}
\put(144,18){$($} \put(155,18){$)$} \put(164,18){$($}
\put(174,18){$($} \put(196,18){$)$}

\end{picture}
\end{center}

\begin{center}
\begin{picture}(200,45)
\put(94,55){$\vector(0,-1){10}$}

\put(1,8){\small 1}    \put(11,8){\small 2} \put(21,8){\small 3}
\put(31,8){\small 4}   \put(41,8){\small 5}  \put(51,8){\small 6}
\put(61,8){\small 7}   \put(71,8){\small 8} \put(81,8){\small 9}
\put(91,8){\small 10}

\put(103,8){\small{-1}} \put(113,8){\small{-2}}
\put(123,8){\small{-3}} \put(133,8){\small{-4}}
\put(143,8){\small{-5}} \put(153,8){\small{-6}}
\put(163,8){\small{-7}} \put(173,8){\small{-8}}
\put(183,8){\small{-9}} \put(193,8){\small{-10}}

\multiput(2,15)(10,0){9}{\circle*{1}} \put(93,15){\circle*{1}}
\multiput(106,15)(10,0){9}{\circle*{1}}\put(198,15){\circle*{1}}

\qbezier[500](2,15)(84,65)(166,15)
\qbezier[500](12,15)(17,24)(22,15)
\qbezier[500](32,15)(84,55)(136,15)
\qbezier[500](42,15)(47,24)(52,15)
\qbezier[500](62,15)(84,40)(106,15)
\qbezier[500](72,15)(77,24)(82,15)
\qbezier[500](82,15)(87.5,24)(93,15)
\qbezier[500](116,15)(121,24)(126,15)
\qbezier[500](146,15)(151,24)(156,15)
\qbezier[500](176,15)(181,24)(186,15)
\qbezier[500](186,15)(192,24)(198,15)

\end{picture}
\end{center}
\caption{The correspondence between a pair $(L,R)$ and a noncrossing
partition.}\label{figLR}
\end{figure}

Conversely, given a type $B$ noncrossing partition, the pair of
subsets $(L, R)$ can be easily determined by the absolute values of
the departure points and the destination points. Thus we have
obtained the desired bijection. An illustration of this
correspondence is given in Figure \ref{figLR}.

%
%
%

\noindent{\it Combinatorial Proof of \eqref{i1}.} We assign the
weight $x$ to departure points and singletons and the weight $1$ to
other points (in the canonical representation). Thus the left-hand
side of  \eqref{i1} equals the total weight of the set
$NC^B(n)$.

To give a combinatorial interpretation of the right-hand side, let
$S_k$ denote the set of partitions in  $NC^B(n)$ with exactly $k$
pairs of nonzero blocks but no singletons, and let $T_k$ denote the
set of partitions in $NC^B(n)$ with exactly $k$ pairs of nonzero
blocks which have at leat two elements. Clearly, every nonzero block
of a partition in $S_k$ contains at least two elements. By
Proposition \ref{pro}, the total weight of $S_k$ equals ${n\choose
2k}{2k\choose k}x^k$. From the above description of the construction
for Proposition \ref{pro}, it is clearly seen that a zero point (in
the canonical representation) can not be covered by an arc of a
nonzero block, or more precisely, a zero point cannot appear between
a pair of departure point and destination point. Furthermore, a
partition in $T_k$ can be obtained by changing some zero points and
transient points to singletons. Conversely, given a partition in
$T_k$, there is only one way to change it back to a partition in
$S_k$ by switching every singleton to either a zero point or a
transient point. This implies that the total weight of $T_k$ equals
the total weight of $S_k$ multiplied by the factor $(1+x)^{n-2k}$.
This completes the proof. \qed

Setting $x=1$, the identity  \eqref{i1} takes the form
\begin{eqnarray}\label{i3}
\sum\limits_{m=0}^{\lfloor \frac{n}{2}\rfloor}{n\choose
2k}{2k\choose k}2^{n-2k}={2n\choose n},
\end{eqnarray}
which can be regarded as a type $B$-analogue of Touchard's formula
\begin{eqnarray*}
\sum\limits_{k=0}^{\lfloor{(n-1)/2}\rfloor}{n-1\choose
2k}C_k2^{n-1-2k}=C_n.
\end{eqnarray*}
 Simion \cite{Simon4} has given a combinatorial interpretation of \eqref{i3} by means  of
a symmetric boolean  decomposition of the lattice $NC^B(n)$  with
the refinement order. Our combinatorial interpretation of (\ref{i2})
will be based on the set $V_n$ of colored type $B$ noncrossing
partitions, which consists of all type $B$ noncrossing partitions of
$[n]$ in which a pair of antipodal singletons may be colored with
two colors, say, black and white.

\noindent{\it Combinatorial Proof of \eqref{i2}.}
 We assign the weight $x^2$ to departure points and singletons,
and assign the weight $(1+x)^2$ to zero points, transient points,
and destination points. In this way, the left-hand side of
\eqref{i2} equals the total weight of the set $NC^B(n)$.

As far as  the right-hand side is concerned, we need to classify the
set $NC^B(n)$ as follows. Given a partition $\pi \in NC^B(n)$, let
$L_\pi$ and $R_\pi$ denote the sets of departure points and
destination points respectively in the canonical representation of
$\pi$. Two partitions $\pi$ and $\sigma$ are in the same class if
$(L_{\pi}, R_{\pi})=(L_{\sigma}, R_{\sigma})$.  Suppose that $(L,R)$
is a pair of feasible sets of departure points and destination
points, namely, there exists $\pi$ such that $(L, R)=(L_\pi,
R_{\pi})$. Let $G(L,R)$ be the set of partitions $\pi\in NC^B(n)$
such that $(L_{\pi}, R_{\pi})=(L,R)$. We further assume that  both
$L$ and $R$ contain $k$ elements.
 Since the departure points and destination points always appear
 in pairs and $x^2(1+x)^2=(x(1+x))^2$, the above weight
assignment is equivalent to the effect of assigning the weight
$x(1+x)$ to  both departure points and destination points. By the
same argument in the proof \eqref{i1}, the noncrossing property
implies that the other points can be either singletons or
non-singletons (transient points or zero points), we deduce that the
total weight of $G(L, R)$ equals
\begin{equation}\label{gk}
(x(1+x))^{2k}(x^2+(1+x)^2)^{n-2k}.
\end{equation}

We next proceed to give a alternative interpretation of the total
weight (\ref{gk}) in terms of colored type $B$ noncrossing
partitions in $V_n$. Assign the weight $x(1+x)$ to black singletons,
zero points, transient points, departure points, destination points
and the weight $1$ to white singletons.  Let us define  $H(L,R)$ for
$V_n$ in the same way as we have defined  $G(L,R)$, that is, the set
of colored type $B$ noncrossing partitions $\pi$ such that $(L_\pi,
R_\pi)=(L,R)$ and $|L|=|R|=k$. By the same argument we see that the
total weight of $H(L,R)$ equals
\begin{equation}\label{hk}
(x(1+x))^{2k}(1+x(1+x)+x(1+x))^{n-2k},
\end{equation}
 which apparently coincides with (\ref{gk}).
This implies that the right-hand side of (\ref{i2}) can be
reinterpreted in terms of partitions in $V_n$.
 It is necessary  to show that the total
weight of $V_n$ equals the right-hand side of  \eqref{i2}. To
construct a partition of $V_n$ with $n-k$ white singletons, we may
first choose $n-k$ white singletons from $[n]$ in ${n \choose k}$
ways. Observe that white singletons have weight $1$ and  the
remaining $2k$ elements can form a type $B$ noncrossing partition
with each element (in the canonical representation) having weight
$x(1+x)$. Clearly, there are ${2k\choose k}$ choices of such
partitions and the weight of each  equals $x^k(1+x)^k$. This
completes the proof.\qed

When $x=-1$, the identity \eqref{i1} becomes
\begin{equation}\label{i4}
\sum\limits_{k=0}^{n}{(-1)^k{n\choose k}^2}=\left\{
  \begin{array}{ll}
  (-1)^r{2r\choose r}, &\mbox{if}\,n=2r;\\[5pt]
  0, &\mbox{otherwise}.
  \end{array}
  \right.
\end{equation}

 We shall provide an involution for the identity \eqref{i4}. It
 should be noted that \eqref{i4} is a type $B$ analogue of the
 identity on the alternating sum of Narayana numbers,
\begin{equation}\label{sign na}
\sum\limits_{k=0}^{n}(-1)^k\frac{1}{n}{n\choose k-1} {n\choose k}
=\left\{
  \begin{array}{ll}
   (-1)^{r+1}C_r, &\mbox{if}\,n=2r+1;\\[5pt]
    0, &\mbox{otherwise}.
  \end{array}
  \right.
\end{equation}

The above identity \eqref{sign na} was first discovered by Bonin,
Shapiro and Simion \cite{Bon} in their study of Schr{\"o}der paths.
It also has been studied by Coker \cite{Coker},  Klazar \cite{Kla},
Eu, Liu and Yeh \cite{Eu}. A combinatorial proof of \eqref{sign na}
has been given by Chen, Shapiro and Yang \cite{Chen1} by using plane
trees and $2$-Motzkin paths.

In the case of $x=-\frac{1}{2}$,  \eqref{i2} reduces to
\begin{eqnarray}\label{twice}
\sum\limits_{k=0}^{n}(-1)^k{n\choose k}{2k\choose
k}4^{n-k}={2n\choose n}.
\end{eqnarray}
This identity can be found in Riordan \cite{Rio} and it was derived
by means of generating functions.

We now  give a combinatorial interpretation of \eqref{i4}. Let
$NC^B_e(n)$ (resp. $NC^B_o(n)$) denote the number of type $B$
noncrossing partitions of $[n]$ into even (resp. odd) pairs of
nonzero blocks. We give a parity reversing involution on type $B$
noncrossing partitions which implies the following formulation of
\eqref{i4}:
\begin{equation}\label{i5}
NC^B_e(n)-NC^B_o(n)=\left\{
  \begin{array}{ll}
  (-1)^r{2r\choose r}, &\mbox{if}\,n=2r;\\[5pt]
  0, &\mbox{otherwise}.
  \end{array}
  \right.
\end{equation}

Let $A_n$ denote the set of type $B$ noncrossing partitions of $[n]$
without the zero block such that every nonzero block contains
exactly two elements. Notice that $A_n$ is empty if $n$ is odd and
$$|A_{2n}|={2n\choose n}.$$ Define the parity of a type
$B$ noncrossing partition as the parity of the number of pairs of
nonzero blocks. Moreover, we define the sign of a partition as $-1$
if it is odd, and as $1$ if it is even.
\begin{thm} \label{rho}
There is a parity reversing involution $\rho$ on the set $NC^B(n)
\setminus A_n$.
\end{thm}

\pf Let $\pi=B_0/B_1/\cdots/B_k\in NC^B(n) \setminus A_n$ be the
canonical representation of  $\pi$. We first list the elements in
$B_0, B_1, \ldots, B_k$ in the increasing order of their absolute
values. Then we define the critical point $i$ of $\pi$ as the first
element in the above order that  is neither a departure point nor a
destination point, or equivalently, is either a zero point,  or a
transient point, or a singleton. We now conduct the
 following map $\rho$ based on two cases concerning the critical
 point:

If the critical point $i$ is a zero point or a transient point, we
take the element $i$ out of  the block  and form a singleton block
$\{i\}$.

If the critical point $i$ is a singleton, we need to  determine
whether $i$ should be put into the zero block or a nonzero block. We
may pay attention to the arc between $i$ and $-i$ in the linear
representation of $\pi$. If this arc does not cross any nonzero
block, then we put $i$ into the zero block. Otherwise, put $i$ into
the nonzero block which has an arc covering $i$ and there are no
arcs of other blocks covered  by this arc.

It is not hard to see that the above map  changes the number of
blocks by one. Moreover, the critical point remains unchanged under
the above map. Thus we can infer that $\rho$ is a parity reversing
involution.
 \qed

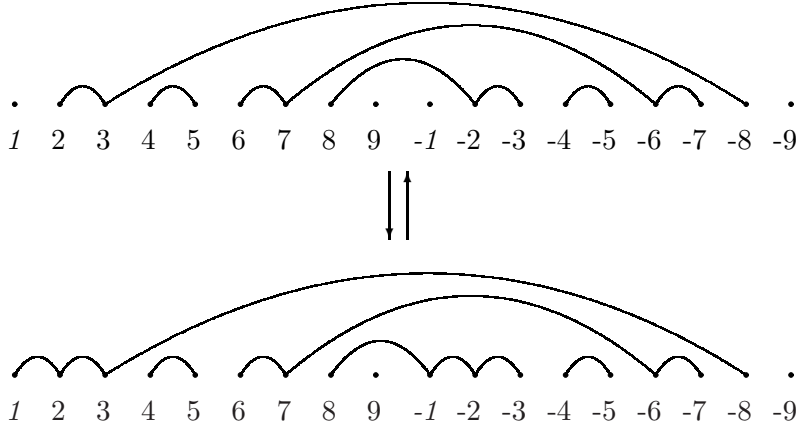
\begin{figure}[h,t]
\setlength{\unitlength}{0.6mm}
\begin{center}
\begin{picture}(180,100)

\put(0,0){\emph{\small1}}

\put(10,0){\small 2}    \put(20,0){\small 3}   \put(30,0){\small 4}
\put(40,0){\small 5}    \put(50,0){\small 6}   \put(60,0){\small 7}
\put(70,0){\small 8}    \put(80,0){\small 9}

\put(90,0){\emph{\small-1}}

\put(100,0){\small{-2}}  \put(110,0){\small{-3}}
\put(120,0){\small{-4}}  \put(130,0){\small{-5}}
\put(140,0){\small{-6}}  \put(150,0){\small{-7}}
\put(160,0){\small{-8}}  \put(170,0){\small{-9}}
\multiput(2,10)(10,0){9}{\circle*{1}}
\multiput(94,10)(10,0){9}{\circle*{1}}
\qbezier[1000](2,10)(7,18)(12,10)
\qbezier[1000](12,10)(17,18)(22,10)
\qbezier[1000](22,10)(94,55)(164,10)
\qbezier[1000](32,10)(37,18)(42,10)
\qbezier[1000](52,10)(57,18)(62,10)
\qbezier[1000](62,10)(103,45)(144,10)
\qbezier[1000](72,10)(83,25)(94,10)
\qbezier[1000](94,10)(99,18)(104,10)
\qbezier[1000](104,10)(109,18)(114,10)
\qbezier[1000](124,10)(129,18)(134,10)
\qbezier[1000](144,10)(149,18)(154,10)
\put(85,55){\vector(0,-1){15}} \put(89,40){\vector(0,1){15}}

\put(0,60){\emph{\small1}}

\put(10,60){\small 2}    \put(20,60){\small 3} \put(30,60){\small 4}
\put(40,60){\small 5}    \put(50,60){\small 6} \put(60,60){\small 7}
\put(70,60){\small 8}    \put(80,60){\small 9}

\put(90,60){\emph{\small-1}}

\put(100,60){\small-2} \put(110,60){\small-3} \put(120,60){\small-4}
\put(130,60){\small-5} \put(140,60){\small-6} \put(150,60){\small-7}
\put(160,60){\small-8} \put(170,60){\small-9}
\multiput(2,70)(10,0){9}{\circle*{1}}
\multiput(94,70)(10,0){9}{\circle*{1}}
\qbezier[1000](12,70)(17,78)(22,70)
\qbezier[1000](22,70)(94,115)(164,70)
\qbezier[1000](32,70)(37,78)(42,70)
\qbezier[1000](52,70)(57,78)(62,70)
\qbezier[1000](62,70)(103,105)(144,70)
\qbezier[1000](72,70)(88,90)(104,70)
\qbezier[1000](104,70)(109,78)(114,70)
\qbezier[1000](124,70)(129,78)(134,70)
\qbezier[1000](144,70)(149,78)(154,70)

\end{picture}
\end{center}
\caption{Involution $\rho$ on a type $B$ noncrossing
partition.}\label{fig6}
\end{figure}

Figure \ref{fig6} gives an example of the involution $\rho$, where
the italic $1$ is the critical point.

Since $A_{2n+1}$ is empty and any partition in $A_{2n}$ has $n$
pairs of nonzero blocks,  the involution $\rho$ leads to a
combinatorial proof of \eqref{i5}.

Appealing to the correspondence between plane trees and $2$-Motzkin
paths, Chen, Deutsch and Elizalde \cite{Chen2} obtained the
following refinements of Coker's identities:
\begin{equation}\label{r1}
\sum\limits_{i=1}^{n}\sum\limits_{j=0}^{n-2i+1}\frac{1}{n}{n\choose
i}{n-1\choose j}{n-i-j\choose
i-1}x^{i-1}y^{j}=\sum\limits_{k=0}^{\lfloor
(n-1)/2\rfloor}C_k{n-1\choose 2k}x^k(1+y)^{n-2k-1},
\end{equation}
\begin{equation}\label{r2}
\sum\limits_{i=1}^{n}\sum\limits_{j=0}^{n-2i+1}\frac{1}{n}{n\choose
i}{n-1\choose j}{n-i-j\choose i-1}x^{2(i-1)}y^{j}z^{n-2i-j+1}
=\sum\limits_{k=0}^{n-1}C_{k+1} {n-1\choose k}x^k(y+z-2x)^{n-1-k}.
\end{equation}

The following  refinements of \eqref{i1} and \eqref{i2} can be
treated as type $B$ analogues of (\ref{r1}) and (\ref{r2}).

\begin{thm}
For $n \geq 1$, we have
\begin{equation}\label{re1}
\sum\limits_{i=0}^{n}\sum_{j=0}^{\lfloor \frac{n-i}{2} \rfloor}
{n\choose i}{n-i \choose j}{n-j-i \choose j}x^jy^i=
\sum\limits_{k=0}^{\lfloor \frac{n}{2}\rfloor}{n\choose
2k}{2k\choose k}x^k(1+y)^{n-2k},
\end{equation}
\begin{equation}\label{re2}
\sum\limits_{i=0}^{n}\sum_{j=0}^{\lfloor \frac{n-i}{2} \rfloor}
{n\choose i}{n-i \choose j}{n-j-i \choose j}x^{2j}y^iz^{n-2j-i}=
\sum\limits_{k=0}^{n}{n\choose k}{2k\choose k}x^k(y+z-2x)^{n-k}.
\end{equation}
\end{thm}

\noindent{\it Proof.} We first consider (\ref{re1}) and shall also
use  weighted type $B$ noncrossing partitions to give the
combinatorial interpretations of both sides of \eqref{re1}. The
weight assignment is almost the same as in the proof of identity
\eqref{i1} except that  a singleton is endowed with the weight $y$.
Suppose that $\pi$ is a type $B$ noncrossing partition with $i$
singletons (in the canonical representation, to be precise). The
remaining $n-i$ elements form a type $B$ noncrossing partition
$\sigma$ such that each nonzero block contains at leat two elements
(in the canonical representation as well). Assume that $\sigma$
contains $j$ nonzero blocks (in the canonical representation). By
Proposition \ref{pro}, there are ${n-i \choose j}{n-j-i \choose
j}={n-i\choose 2j}{2j \choose j}$ choices for $\sigma$.
 According to the weight assignment, $\sigma$ has  weight $x^j$. In view
of the number of singletons, we see that the left-hand side of
 \eqref{re1} equals the sum of the weights of all type $B$
noncrossing partitions on $[n]$.

Regarding the right-hand side, let $S_k$ be the set of the
partitions in $NC^B(n)$ with $k$  pairs of nonzero blocks but no
singletons, and let $T_k$ be the set of partitions in $NC^B(n)$ with
exactly $k$ pairs of nonzero and nonsingleton blocks. The total
weight of any partition in $S_k$ a partition is $x^k$. Meanwhile,
there are ${n\choose 2k}{2k\choose k}$ partitions in $S_k$. Using
the same argument as in the proof of (\ref{i1}), we see that the
total weight of $T_k$ equals the total weight of $S_k$ multiplied by
the factor $(1+y)^{n-2k}$. Thus the right-hand side of \eqref{re1}
also equals the total weight of $NC^B(n)$. Hence (\ref{re1}) is
proved.

We now turn to the proof of  \eqref{re2} using a different weight
assignment. Assign the weight $x$ to departure points and
destination points, the weight $y$ to singletons and the weight $z$
to zero points, transient points. Suppose that $\pi$ is a partition
with $i$ singletons and $j$ nonzero blocks with at least two
elements. Since the departure points and destination points appear
in pairs, there are $n-2j-i$ other elements with weight $z$. From
Proposition \ref{pro} it follows that the left-hand side of
\eqref{re2} equals the  total weight of the set $NC^B(n)$.

On the other hand, we may consider the set $V_n$ of colored type $B$
noncrossing partitions, as defined before. Assign the weight $x$ to
black singletons, zero points, transient points, departure points,
destination points and the weight $y+z-2x$ to white singletons. Let
$G(L,R)$ and $H(L,R)$ denote the subsets of $NC^B(n)$ and $V_n$
respectively as in the proof of the identity \eqref{i2}. According
to the weight assignment, the departure points and destination
points have weight $x$, singletons have weight $y$ and the other
points have weight $z$, thus the total weight of the set $G(L,R)$ is
$x^{2k}(y+z)^{n-2k}$. Similarly, the total weight of $H(L,R)$ equals
$x^{2k}(x+x+y+z-2x)^{n-2k}$, which coincides with the total weight
of  $G(L,R)$. Now it suffices to show that the total weight of $V_n$
agrees with  the right-hand side of  \eqref{i2}. In order to compute
the total weight of $V_n$, we can first choose $n-k$ white
singletons from $[n]$ in ${n \choose k}$ ways. These white
singletons have total weight  $(y+z-2x)^{n-k}$. The remaining $2k$
elements form a type $B$ noncrossing partition with each kind of
points having weight $x$. There are  ${2k\choose k}$ choices of such
partitions and the weight of each partition equals $x^k$.
 This completes the proof.\qed

\section{ Type $B$ Analogues of the Identities of Mansour and Sun }\label{Daw}

In this section, we provide   type $B$ analogues  along with
combinatorial proofs of the identities \eqref{m1}, \eqref{m2} and
\eqref{m3} due to Mansour and Sun \cite{Man}. We begin with the
following identity which can be considered as a type $B$ analogue of
\eqref{m1}.

\begin{thm} \label{bm1}
For $n\geq 0$, we have
\begin{equation}\label{i6}
\sum\limits_{k=0}^n (-1)^{n-k}{n \choose k}P_{k}(x)(1+x)^{n-k}=
\left\{
  \begin{array}{cl}
  x^r{2r\choose r}, &\mbox{if }\,n=2r;\\[6pt]
  0, &\mbox{otherwise}.
  \end{array}
  \right.
\end{equation}
\end{thm}

Setting $x=1$, \eqref{i6} reduces to the following identity of
Dawson, see, Riordan \cite[p.71]{Rio},
\begin{equation}\label{Reed}
\sum\limits_{k=0}^n (-1)^{n-k}{n \choose k}{2k\choose k}2^{n-k}=
\left\{
  \begin{array}{cl}
  {2r\choose r}, &\mbox{if}\,n=2r;\\[6pt]
  0, &\mbox{otherwise}.
  \end{array}
  \right.
\end{equation}
 Andrews \cite{And} has given a proof of (\ref{Reed}) by employing
  Gauss's second summation theorem, which is stated as
\[_2F_1\left[\begin{array}{c}
  a, b\\
  1/2+a/2+b/2\\
  \end{array}; 1/2\right]=
  \frac{\Gamma({1/ 2})\Gamma(1/2+a/2+b/2)}
  {\Gamma(1/2+a/2)\Gamma(1/2+b/2)}.\]

 \noindent{\it
Combinatorial Proof of Theorem \ref{bm1}.}  The left-hand side of
(\ref{i6}) can be interpreted in terms of weighted type $B$ colored
partitions in $V_n$. We use the following weight assignment.
Transient points, zero points and destination points are given the
weight $1$; black singletons and departure points are given the
weight $x$;  white singletons are given  the weight $-1-x$, which is
equivalent to the assignment of either $-x$ or $-1$.

To construct a partition in  $V_n$,  we first choose a subset $S$ of
$n-k$ elements from $[n]$  to form  white singletons with weight
$-1-x$. There are ${n \choose k}$ ways to choose $S$, and these
white singletons will contribute a factor $(-1-x)^{n-k}$ to the
 weight. On the other hand,  the remaining elements constitute
 a noncrossing partition on $2k$ elements with the same
weight assignment except that the singletons are assumed to be
black. It follows that such partitions have a total weight $P_k(x)$.
Thus the left-hand side of \eqref{i6} equals the total weight of the
set $V_n$.

We now aim to construct a sign reversing involution $\theta$ on $V_n
\setminus A_n$, where $A_n$ is defined in the preceding section. By
the definition of $A_n$, we may regard the partitions in $A_n$ as
noncrossing partitions in $V_n$ with no zero block, no singletons
and no transient points. The involution $\theta$ can be described as
follows.

Given a type $B$ colored noncrossing partition $\pi\in V_n \setminus
A_n$, in the increasing order of the absolute values, we search the
first point $i$ in the canonical representation of $\pi$ that is
neither a departure point nor a destination point of $\pi$. Then we
conduct the following operations:
\begin{itemize}
\item
If $i$ is a black singleton, then we change it to a white singleton
with weight $-x$;

\item
If $i$ is a transient point or a zero point, then we change it to a
white singleton with weight $-1$;

\item
If $i$ is a white singleton with weight $-x$, then we change it to a
black singleton;

\item
If $i$ is a white singleton with weight $-1$, then we change it to a
transient point or a zero point by  the same criterion in the proof
of Theorem \ref{rho}.
\end{itemize}

Evidently, $\theta$ is weight preserving and sign reversing. It is
also easily seen the critical point after the map is neither a
departure point nor a destination point. Moreover, in the increasing
order of the absolutes, the elements before $i$ stay unchanged since
the above operations do not cause additional departure point or
destination point. Therefore, the critical point remains the same.
Thus $\theta$ is a sign reversing and weight preserving involution.

Our final task is to compute the total weight of $A_n$. First
consider the case when $n$ is odd. For a type $B$ noncrossing
partition $\pi$ on  $[n]$, it is clear that there exists at least
one point which is neither a departure point nor a destination
point. By the involution $\theta$  we deduce that the sum of weights
of type $B$ colored noncrossing partitions on $[n]$ equals zero. For
$n=2r$, the total weight of $V_n$ equals to the total weight of
$A_n$, which equals $x^r{2r\choose r}$. This completes the
proof.\qed

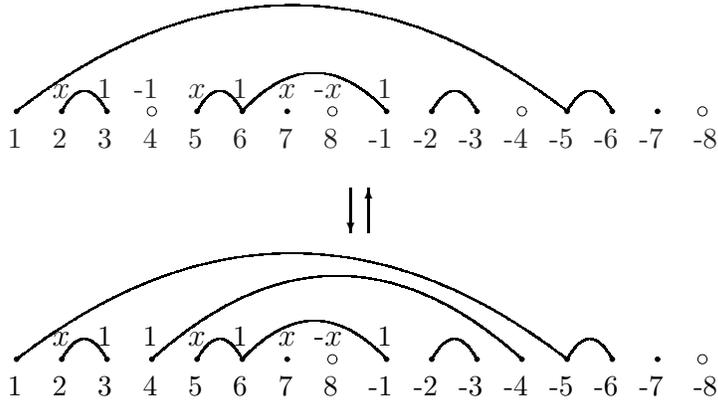
\begin{figure}[h,t]
\setlength{\unitlength}{0.6mm}
\begin{center}
\begin{picture}(150,80)
\put(0,47){\small 1}    \put(10,47){\small2}   \put(20,47){\small3}
\put(30,47){\small4}   \put(40,47){\small5}   \put(50,47){\small6}
\put(60,47){\small7}   \put(70,47){\small8}   \put(80,47){\small-1}
\put(90,47){\small-2}  \put(100,47){\small-3} \put(110,47){\small-4}
\put(120,47){\small-5} \put(130,47){\small-6} \put(140,47){\small-7}
\put(152,47){\small-8}

\put(10,58){$x$}  \put(20,58){\small1}   \put(28,58){\small-1}
\put(40,58){$x$} \put(50,58){\small1}   \put(60,58){$x$}
\put(68,58){-$x$} \put(82,58){\small1}

\multiput(2,55)(10,0){3}{\circle*{1}} \put(32,55){\circle{2}}
\multiput(42,55)(10,0){3}{\circle*{1}} \put(72,55){\circle{2}}

\multiput(84,55)(10,0){3}{\circle*{1}} \put(114,55){\circle{2}}
\multiput(124,55)(10,0){3}{\circle*{1}} \put(154,55){\circle{2}}

 \qbezier[1000](2,55)(63,102)(124,55)
\qbezier[1000](12,55)(17,64)(22,55)
\qbezier[1000](42,55)(47,64)(52,55)
\qbezier[1000](52,55)(68,72)(84,55)

\qbezier[1000](94,55)(99,64)(104,55)
\qbezier[1000](124,55)(129,64)(134,55)

\put(76,38){\vector(0,-1){10}}\put(80,28){\vector(0,1){10}}

\put(0,-8){\small1}    \put(10,-8){\small2}   \put(20,-8){\small3}
\put(30,-8){\small4}   \put(40,-8){\small5}   \put(50,-8){\small6}
\put(60,-8){\small7}   \put(70,-8){\small8}   \put(80,-8){\small-1}
\put(90,-8){\small-2}  \put(100,-8){\small-3} \put(110,-8){\small-4}
\put(120,-8){\small-5} \put(130,-8){\small-6} \put(140,-8){\small-7}
\put(152,-8){\small-8} \multiput(2,0)(10,0){7}{\circle*{1}}
\put(72,0){\circle{2}} \multiput(84,0)(10,0){3}{\circle*{1}}
\put(114,0){\circle*{1}} \multiput(124,0)(10,0){3}{\circle*{1}}
\put(154,0){\circle{2}}

\qbezier[1000](2,0)(63,47)(124,0) \qbezier[1000](12,0)(17,9)(22,0)
\qbezier[1000](32,0)(73,37)(114,0) \qbezier[1000](42,0)(47,9)(52,0)
\qbezier[1000](52,0)(68,17)(84,0)

\qbezier[1000](94,0)(99,9)(104,0)
\qbezier[1000](124,0)(129,9)(134,0)

\put(10,3){$x$}  \put(20,3){\small1}   \put(30,3){\small1}
\put(40,3){$x$}  \put(50,3){\small1}   \put(60,3){$x$}
\put(68,3){-$x$} \put(82,3){\small1}

\end{picture}
\end{center}
\caption{The involution $\theta$ on a type $B$ colored noncrossing
partition.} \label{fig7}
\end{figure}

Figure \ref{fig7} is an illustration of the involution $\theta$,
where a circle stands for a white singleton.

We now turn to a type $B$ analogue of the identity \eqref{m2}.

\begin{thm} For $n\geq 0$, we have
\begin{eqnarray}\label{i7}
\sum_{k=0}^n (-1)^{n-k}{n \choose
k}P_{k}(x^2)(1-x)^{2(n-k)}=x^{n}{2n \choose n}.
\end{eqnarray}
\end{thm}
It should be noted that setting $x=-1$ in \eqref{i7}, we arrive at
\eqref{twice} again.

\noindent{\it Combinatorial Proof of \eqref{i7}.} The proof is
similar to that  of \eqref{i6}. We still consider  type $B$ colored
noncrossing partitions with  a different weight assignment.
Destination points, zero points and  transient points are given the
weight $1$; departure points and  black singletons are given the
weight $x^2$;  white singletons are given the weight $-1+2x-x^2$,
or, equivalently,  the weight of a white singleton can be either
$-1$, or $2x$ or $-x^2$. Then the left-hand side of \eqref{i7} can
be interpreted as  the total weight of the set $V_n$.

Let ${D}_n$ be the set of colored noncrossing partitions in $V_n$
that have only two types of points (in the canonical
representation): (1) departure points or destination points. (2)
white singletons with weight $2x$. In other words, in any partition
in ${D}_n$, the following four types of points are not allowed: (1)
zero points; (2) transient points; (3) black singletons; (4) white
singletons with weight $-1$ or $-x^2$.  To prove \eqref{i7}, we
proceed to construct  a sign reversing and weight preserving
involution $\eta$ on the set $V_n \setminus {D}_n$.

 Given a type $B$ colored noncrossing partition $\pi\in V_n
 \setminus
{D}_n$, we seek the first point $i$ in the increasing order of the
absolute values which is neither a departure point nor a destination
point. As usual, $i$ is called the critical point. The map $\eta$ is
defined by the following operations:
\begin{itemize}
\item
If $i$ is a black singleton, then  set $i$ to be a white singleton
with weight $-x^2$;

\item
If $i$ is a transient point or a zero point, then set $i$ to be a
white singleton with weight $-1$;

\item
If $i$ is a white singleton with weight $-x^2$, then set $i$  to be
a black singleton;

\item
If $i$ is a white singleton with weight $-1$, then set $i$ to be a
transient point or a zero point according to the criterion as given
before.
\end{itemize}
It can be verified that $\eta$ is a sign reversing and weight
preserving involution. Thus the total weight of $V_n$ equals the
total weight of ${D}_n$.

We now are left the task  to show that the total weight of the set
${D}_n$ equals the right-hand side  of (\ref{i7}), namely, $x^{n}{2n
\choose n}$.
 To construct a weighted type $B$ colored noncrossing partition in
${D}_n$, we can first choose $2k$ elements from $[n]$ to construct a
type $B$  noncrossing partition with $k$ departure points and $k$
destination points. There are ${n \choose 2k}{2k \choose k}$
choices.  The remaining $n-2k$ points are taken to be white
singletons with weight $2x$. Thus the total weight of the set
${D}_n$ equals
\[ \sum_{k=0}^{\lfloor n/2 \rfloor} {n \choose 2k}{2k \choose k}x^{2k}(2x)^{n-2k},\]
which can be rewritten as
\[x^n \sum\limits_{k=0}^{\lfloor n/2 \rfloor} {n \choose 2k}{2k \choose k} 2^{n-2k}.\]
Invoking  the  identity \eqref{i3} gives \eqref{i7}, so the proof is
complete. \qed

We remark in passing that \eqref{m3} is related to the following
identity obtained by  Chen and Pang \cite{Chen4}
\begin{equation}\label{cp}
\sum\limits_{k=0}^{n}\frac{1}{n}{n\choose k-1}{n\choose
k}x^k(x+1)^{n-k} =\sum\limits_{k=0}^{n}(-1)^{n-k}{n+k\choose
{n-k}}\frac{1}{k+1}{2k\choose k}(x+1)^{k},
\end{equation}
which was independently derived by Mansour and Sun \cite{MS} in a
slightly different form
\begin{equation}\label{sun}
\sum\limits_{k=0}^{n}\frac{1}{n}{n\choose k-1}{n\choose k}x^k
=\sum\limits_{k=0}^{n}{n+k\choose {n-k}}\frac{1}{k+1}{2k\choose
k}(x-1)^{n-k}.
\end{equation}
Upon substituting $x$ with $x/(1-x)$, we see that (\ref{cp}) can be
recast in the form of (\ref{sun}), that is,
\begin{equation} \label{ssun}
\frac{N_n(x)}{(x-1)^n}=\sum\limits_{k=0}^{n}{n+k\choose{n-k}}
C_k\frac{1}{(x-1)^k}.
\end{equation}
Then the identity (\ref{m3})  follows from the Legendre inversion
formula \cite{Rio}
\[ a_n=\sum_{k=0}^n {n+k \choose n-k}b_k   ~~\Longleftrightarrow~~
b_n=\sum_{k=0}^n (-1)^{n-k}\frac{2k+1}{2n+1}{2n+1 \choose n-k}a_k.\]
 Below is a type $B$ analogue of (\ref{cp}), which
yields  a type $B$ analogue of (\ref{m3}).

\begin{thm}  For $n\geq 0$, we have
\begin{equation}\label{bi}
\sum\limits_{k=0}^{n}{n\choose k}^2x^k(x-1)^{n-k}=
\sum\limits_{k=0}^{n}(-1)^{n-k}{n+k\choose k}{n\choose k}x^k.
\end{equation}
\end{thm}

\noindent{\it Proof.} The weight assignment of the points of a type
$B$ noncrossing partition  on $[n]$ is given as follows. A departure
point or a singleton is endowed with the weight $x$, whereas the
other kinds of points are given  the weight $x-1$. In this way,
 the left-hand side of (\ref{bi}) equals the total weight of the
set $NC^B(n)$.

We need an equivalent weight assignment, that is,  departure points
and singletons always have weight $x$ and the weight $x-1$ for other
types of points can be interpreted as the assignment of either $x$
or $-1$. Using this weight assignment, we proceed to show that the
summand of the right-hand side equals the sum of weights of type $B$
noncrossing partitions on $[n]$ with exactly $k$ elements having
weight $x$. Such partitions can be constructed as follows. We first
choose a $k$-subset $S$ of $[n]$ under the condition that if an
element $i$ or $-i$ has weight $x$, then $i\in S$. Then we choose
$m$ elements from $S$ to be departure points or singletons, denoted
by $L=\{l_1, l_2, \ldots, l_m\}$. Meanwhile,  we choose another set
of $m$ elements from the set $[n]$, denoted by $R=\{r_1, r_2,
\ldots, r_m\}$. Applying the bijection of Reiner \cite{Reiner} to
the pair $(L,R)$, we get a noncrossing partition with $m$ pairs of
antipodal nonzero blocks.

It remains to compute the sum of weights of such partitions. Note
that the departure points and singletons always have weight $x$ and
consequently belong to $S$. It suffices to determine  which of the
other types of points can be given weight $x$. The weight of the
destination points, zero points and transient points can be either
$x$ or $-1$ subject to the condition for the choice of $S$. If the
absolute value of such an element belongs to the set $S$, then it
has  weight $x$; otherwise, it has weight $-1$. That is to say that
subject to the condition on the choice of $S$, the weight assignment
for all the elements are uniquely determined. So the weight of a
partition constructed by the above procedure equals $(-1)^{n-k}x^k$.
On the other hand, the number of such partitions equals
\[{n\choose k}\sum\limits_{m=0}^{k}{k\choose m}{n\choose m}
 ={n\choose k}{n+k\choose k}.
\]
 This completes the proof.  \qed

In virtue of the identity \eqref{bi}, we can express the central
binomial coefficients in terms of  the Narayana polynomials of type
$B$.

\begin{thm} For $n\geq 0$,
\begin{equation}\label{bi3}
{2n \choose n}=\sum_{k=0}^n\frac{2k+1}{2n+1}{2n+1 \choose
n-k}P_k(x)(1-x)^{n-k}.
\end{equation}
\end{thm}
\pf Substituting $x$ with  $\frac{1}{1-x}$ into \eqref{bi}, and
observing that
\[ {n+k\choose k}{n\choose k}={n+k\choose n-k}{2k\choose k},\]
we find
\begin{equation} \label{bi2}
\sum\limits_{k=0}^{n}{n\choose k}^2x^k=
\sum\limits_{k=0}^{n}(-1)^{n-k}{n+k\choose n-k}{2k\choose
k}(1-x)^{n-k},
\end{equation}
which serves as a type $B$ analogue of (\ref{ssun}). Rewriting
\eqref{bi2} as
\[  \frac{P_n(x)}{(x-1)^n}=  \sum\limits_{k=0}^{n}{n+k\choose n-k}
{2k\choose k}\frac{1}{(x-1)^{k}},\]  applying the Legendre inversion
formula, we arrive at the expression \eqref{bi3}.  \qed

\vspace{.2cm} \noindent{\bf Acknowledgments.} This work was
supported by the 973 Project, the PCSIRT Project of the Ministry of
Education, and the National Science Foundation of China.

\end{document}